\documentclass[11pt,twoside]{amsart}

\usepackage{amsfonts}

\title{A Noether-Deuring theorem for derived categories}

\author{Alexander Zimmermann}
\address{\newline
Universit\'e de Picardie,
\newline D\'epartement de Math\'ematiques et LAMFA (UMR 7352 du CNRS),
\newline 33 rue St Leu,
\newline F-80039 Amiens Cedex 1,
\newline France}
\email{alexander.zimmermann@u-picardie.fr}

\newtheorem*{Theo2}{{Theorem}}
\newtheorem{Lemma1}{{Lemma}}
\newtheorem{Theo1}[Lemma1]{{Theorem}}
\newtheorem{Def1}[Lemma1]{{Definition}}
\newtheorem{Prop1}[Lemma1]{{Proposition}}
\newtheorem{Claim1}[Lemma1]{{Claim}}
\newtheorem{Rem1}[Lemma1]{{Remark}}
\newtheorem{Cor1}[Lemma1]{{Corollary}}
\newtheorem{Ex1}[Lemma1]{{Example}}

\newenvironment{Lemma}{\begin{Lemma1}}{\end{Lemma1}}

\newenvironment{Prop}{\begin{Prop1}}{\end{Prop1}}
\newenvironment{Rem}{\begin{Rem1}\rm}{\end{Rem1}}
\newenvironment{Theorem}{\begin{Theo1}}{\end{Theo1}}

\newenvironment{Cor}{\begin{Cor1}}{\end{Cor1}}

\newcommand{\uar}{\uparrow}

\newcommand{\lra}{\longrightarrow}
\newcommand{\lla}{\longleftarrow}
\newcommand{\ra}{\rightarrow}
\newcommand{\sdp}{\times\kern-.2em\vrule height1.1ex depth-.05ex}
\newcommand{\epi}{\lra \kern-.8em\ra}

\newcommand{\N}{{\mathbb N}}

\newcommand{\dickebox}{{\vrule height5pt width5pt depth0pt}}

 \normalbaselines

\date{November 11, 2011; revised December 28, 2011}

\begin{document}

\begin{abstract}
We prove a Noether-Deuring theorem for the derived category of
bounded complexes of modules over a Noetherian algebra.
\end{abstract}

\maketitle

\section*{Introduction}

The classical Noether-Deuring theorem states that given an algebra $A$
over a field $K$ and a finite extension field $L$ of $K$,  two $A$-modules
$M$ and $N$ are isomorphic as $A$-modules, if $L\otimes_KM$
is isomorphic to $L\otimes_KN$
as an $L\otimes_KA$-module. In 1972 Roggenkamp gave a nice
extension of this result to
extensions $S$ of local commutative Noetherian rings $R$ and
modules over Noetherian $R$-algebras.

For the derived category of $A$-modules no such generalisation
was documented before.
The purpose of this note is to give a version of the
Noether-Deuring theorem, in the
generalised version given by Roggenkamp, for right bounded
derived categories of
$A$-modules. If there is a morphism $\alpha\in Hom_{D(\Lambda)}(X, Y)$,
then it is fairly easy to show that for a faithfully flat ring
extension $S$ over $R$ the fact that $id_S\otimes\alpha$ is an isomorphism
implies that $\alpha$ is an isomorphism. This is done in
Proposition~\ref{alreadydefined}. More delicate is the question if
only an isomorphism in $Hom_{D(S\otimes_R\Lambda)}(S\otimes_RX, S\otimes_RY)$
is given. Then we need further  finiteness conditions on $\Lambda$ and on $R$
and proceed by completion of $R$ and then a classical going-down argument.
This is done in Theorem~\ref{abstractisos} and
Corollary~\ref{finallyfaithfulprojectiveextensions}.

For the notation concerning derived categories we refer
to Verdier~\cite{Verdier}. In particular, $D(A)$
(resp $D^-(A)$, resp $D^b(A)$) denotes the derived category
of complexes (resp. right bounded complexes, resp. bounded complexes)
of finitely generated $A$-modules, $K^-(A-proj)$
(resp. $K^b(A-proj)$, resp $K^{-,b}(A-proj)$)
is the homotopy category of right bounded complexes
(resp. bounded complexes, resp. right bounded complexes with
bounded homology) of finitely generated projective $A$-modules.
For a complex $Z$ we denote by $H_i(Z)$ the homology of $Z$ in degree
$i$, and by $H(Z)$ the graded module given by the homology of $Z$.

\section{The result}

We start with an easy observation.

\begin{Prop}\label{alreadydefined}
Let $R$ be a
commutative ring and let $\Lambda$ be an
$R$-algebra.
Let $S$ be a commutative faithfully flat $R$-algebra.
Denote by $D(\Lambda)$ the derived
category of complexes of finitely generated $\Lambda$-modules.
Then if there is $\alpha\in Hom_{D(\Lambda)}(X,Y)$ so that
$id_{S}\otimes_R^{\mathbb L}\alpha\in
Hom_{D(S\otimes_R\Lambda)}(S\otimes_R^{\mathbb L}X,S\otimes_R^{\mathbb L}Y)$
is an isomorphism in $D(S\otimes_R\Lambda)$,
then $\alpha$ is an isomorphism in $D(\Lambda)$.
\end{Prop}

Proof.
Let $Z$ be a complex in $D(\Lambda)$.
Since $S$ is flat over $R$ the functor $S\otimes_R-:R-Mod\lra S-Mod$
is exact, and hence
the left derived functor $S\otimes_R^{\mathbb L}-$
coincides with the ordinary tensor product functor
$S\otimes_R-$. We can therefore
work with the usual tensor product and a complex $Z$ of $\Lambda$-modules.

We claim that since $S$ is flat, $S\otimes_R-$ induces an isomorphism
$S\otimes_RH(Z)\simeq H(S\otimes_R^{\mathbb L}Z)$.

If $\partial_Z$ is the differential of $Z$, then
$$0\lra\ker(\partial_Z)\lra Z\stackrel{\partial_Z}{\lra}im(\partial_Z)\lra 0$$
is exact in the category of  $\Lambda$-modules.

Since $S$ is flat,
$$0\lra S\otimes_R\ker(\partial_Z)\lra S\otimes_RZ
\stackrel{id_{S}\otimes_R\partial_Z}{\lra}S\otimes_Rim(\partial_Z)\lra 0$$
is exact. Hence
$$\ker(id_{S}\otimes_R\partial_Z)=S\otimes_R\ker(\partial_Z)
\mbox{ and }
im(id_S\otimes_R\partial_Z)=S\otimes_Rim(\partial_Z).$$
This shows the claim.

Since $id_S\otimes_R\alpha$ is an isomorphism, its cone
$C(id_S\otimes_R\alpha)$ is acyclic. Moreover,
$C(id_S\otimes_R\alpha)=S\otimes_RC(\alpha)$ by the very construction
of the mapping cone.
But now,
$$0=H(C(id_S\otimes_R\alpha))=H(S\otimes_RC(\alpha))=S\otimes_RH(C(\alpha)).$$
Since $S$ is faithfully flat, this implies $H(C(\alpha))=0$
and therefore $C(\alpha)$ is acyclic.
We conclude that $\alpha$ is an isomorphism in
$D(\Lambda)$ which shows the statement. \dickebox

\begin{Rem}
Observe that we assumed that $X\stackrel{\alpha}{\lra} Y$ is assumed to
be a morphism in $D(\Lambda)$.
The question if the existence of an
isomorphism $S\otimes_RX\stackrel{\hat\alpha}{\lra} S\otimes_RY$
in $D(S\otimes_R\Lambda)$ implies the existence of
a morphism $\alpha:X\lra Y$ in $D(\Lambda)$
so that $id_S\otimes_R^{\mathbb L}\alpha$ is an isomorphism is left open.
Under stronger hypotheses this is the purpose of
Theorem~\ref{abstractisos} below. The proof follows
\cite{Rog72} which deals with the module case.
\end{Rem}

\begin{Lemma}\label{stoensorhomequalhomtensors}
If $S$ is a faithfully flat $R$-module
and $\Lambda$ is a Noetherian
$R$-algebra, then for all objects $X$ and $Y$ of $D^b(\Lambda)$ we get
$$Hom_{D^b(S\otimes_R\Lambda)}(S\otimes_RX,S\otimes_RY)
\simeq S\otimes_RHom_{D^b(\Lambda)}(X,Y).$$
\end{Lemma}

Proof.
Since $S$ is flat over $R$, the functor
$S\otimes_R-$ preserves quasi-isomorphisms and therefore we get a morphism
$$S\otimes_RHom_{D^b(\Lambda)}(U,V)\lra  Hom_{D^b(S\otimes_R\Lambda)}(S\otimes_RU,S\otimes_RV)$$
in the following way.
Given a morphism $\rho$ in $Hom_{D^b(\Lambda)}(U,V)$ represented by the triple
$\left(U\stackrel{\alpha}{\lla}W\stackrel{\beta}{\lra}V\right)$, for a quasi-isomorphism
$\alpha$ and a morphism of complexes $\beta$, and $s\in S$ then map $s\otimes\rho$ to
$\left(S\otimes_RU
\stackrel{id_S\otimes \alpha}{\lla}S\otimes_RW
\stackrel{s\otimes\beta}{\lra}S\otimes_RV\right)$.
This is natural in $U$ and $V$.

We use the equivalence of categories $K^{-,b}(\Lambda-proj)\simeq D^b(\Lambda)$
and suppose therefore that $X$ and $Y$ are right bounded complexes
of finitely generated projective $\Lambda$-modules.
But
$$S\otimes_RHom_\Lambda(\Lambda^n,U)=S\otimes_RU^n=(S\otimes_RU)^n=
Hom_{S\otimes_R\Lambda}((S\otimes_R\Lambda)^n,S\otimes_RU)$$
which proves the statement in case $X$ or $Y$ is in $K^b(A-proj)$
since then a homomorphism is given by a direct sum of finitely many
homogeneous mappings in those degrees where the complexes do
both have non zero components. Now, tensor product commutes with direct sums.

We come to the general case.
Recall the so-called stupid truncation $\tau_N$ of a complex. Let $Z$
be a complex in $K^{-,b}(\Lambda-proj)$, denote by $\partial$ its differential
and let $N\in\N$ so that $H_n(Z)=0$ for all $n\geq N$.
We denote the homogeneous components of $\partial$ so
that $\partial_n:Z_n\lra Z_{n-1}$ for all $n$.
Let $\tau_NZ$ be the complex given by $(\tau_N Z)_n=Z_n$ if $n\leq N$ and
$(\tau_NZ)_n=0$ else. The differential $\delta$ on $\tau_NZ$ is defined to be
$\delta_n=\partial_n$ if $n\leq N$ and $\delta_n=0$ else.
Now, $\ker(\partial_N)=:C_N(Z)$ is a finitely generated $\Lambda$-module.
Therefore
we get an exact triangle, called in the sequel  the truncation triangle for $Z$,
$$\tau_N Z\lra Z\lra C_N(Z)[N+1]\lra (\tau_N Z)[1]$$
for all objects $Z$ in $K^{-,b}(A-proj)$.
Obviously $\tau_N(S\otimes_RZ)=S\otimes_R\tau_NZ$
and since $S$ is flat over $R$ also $C_N(S\otimes_RZ)=S\otimes_RC_N(Z)$.

We choose $N$ so that $H_n(X)=H_n(Y)=0$ for all $n\geq N$.
To simplify the notation denote for the moment the bifunctor
$Hom_{K^{-,b}(\Lambda-proj)}(-,-)$ by $(-,-)$, the bifunctor $Hom_{K^{-,b}(S\otimes_R\Lambda-proj)}(-,-)$
by $(-,-)_S$ and the bifunctor
$S\otimes_RHom_{K^{-,b}(\Lambda-proj)}(-,-)$ by $S(-,-)$.
Further, put $S\otimes_RX=:X_S$ and $S\otimes_RY=:Y_S$.
From the long exact sequence
obtained by applying $(X_S,-)_S$ to the truncation triangle of $Y_S$,
we get a commutative diagram with exact lines $(\dagger)$
{\footnotesize $$\begin{array}{ccccccccc}
(X_S,C_N(Y_S)[N])_S&\ra& (X_S,\tau_NY_S)_S&\ra&(X_S,Y_S)_S&\ra&
(X_S,C_N(Y_S)[N+1])_S&\ra& (X_S,\tau_NY_S[1])_S\\
\uar&&\uar&&\uar&&\uar&&\uar\\
S(X,C_N(Y)[N])&\ra& S(X,\tau_NY)&\ra&S(X,Y)&\ra&
S(X,C_N(Y)[N+1])&\ra& S(X,\tau_NY[1])
\end{array}$$}
Since $\tau_N(Y_S)$ is a bounded complex of projectives,
$$(X_S,\tau_NY_S)_S=S\otimes_R(X,\tau_NY)
\mbox{ and }(X_S,\tau_NY_S[1])_S=S\otimes_R(X,\tau_NY[1]).$$
We apply $(-,C_N(Y_S)[k])_S$, for a fixed integer $k$, to the truncation triangle for $X_S$
and obtain an exact sequence
\begin{eqnarray*}
(\tau_NX_S[1],C_N(Y_S)[k])_S&\ra& (C_N(X_S)[N+1],C_N(Y_S)[k])_S\ra(X_S,C_N(Y_S)[k])_S
\ra\\
&&\phantom{}
\ra(\tau_NX_S,C_NY_S[k])_S\ra(C_N(X_S)[N],C_N(Y_S)[k])_S
\end{eqnarray*}
and a commutative diagram analogous to the diagram $(\dagger)$.

\medskip

Now, for morphisms between finitely presented
$\Lambda$-modules $M$ and $N$ we do have that the natural map
$$S\otimes_RHom_\Lambda(M,N)\lra Hom_{S\otimes_R\Lambda}(S\otimes_RM,S\otimes_RN)$$
is an isomorphism.  Indeed, let
$$P_1\lra P_0\lra M\lra 0$$
be the first terms of a projective resolution of $M$ as a $\Lambda$-module. Then
$$S\otimes_RP_1\lra S\otimes_RP_0\lra S\otimes_RM\lra 0$$
are the first terms of a projective resolution of $S\otimes_RM$ as an
$S\otimes_R\Lambda$-module, and, denoting by $SM:=S\otimes_RM$, $SN:=S\otimes_RN$, $SP_i:=S\otimes_RP_i$ for $i\in\{0,1\}$, and $S\Lambda:=S\otimes_R\Lambda$, we get
$$
\begin{array}{ccccccc}
Hom_{S\Lambda}(SM,SN)&\hookrightarrow&
Hom_{S\Lambda}(SP_0,SN)&\ra&
Hom_{S\Lambda}(SP_1,SN)\\
\uar&&\uar&&\uar\\
S\otimes_RHom_{\Lambda}(M,N)&\hookrightarrow&
S\otimes_RHom_{\Lambda}(P_0,N)&\ra&
S\otimes_RHom_{\Lambda}(P_1,N)
\end{array}
$$
is a commutative diagram with exact lines. The second and the third vertical
morphisms are isomorphisms. Indeed,
\begin{eqnarray*}
S\otimes_RHom_\Lambda(\Lambda^n,N)&=&S\otimes_RN^n\\
&=&(S\otimes_RN)^n\\
&=&Hom_{S\otimes_R\Lambda}(S\otimes_R\Lambda^n,S\otimes_RN)
\end{eqnarray*}
and since $P_0$ and $P_1$ are direct factors of $\Lambda^n$, for certain $n$
gives the result. Therefore the left most vertical homomorphism is an
isomorphism as well.

Given a projective
resolution $P_\bullet\lra M$ of $M$, denote by
$\partial_n:\Omega^nM\hookrightarrow P_{n-1}$ the embedding of the $n$-th syzygy
of $M$ into the degree $n-1$ homogeneous component of the projective resolution.
Then $$Ext^n_\Lambda(M,N)=
Hom_\Lambda(\Omega^nM,N)/\left(Hom_\Lambda(P_{n-1},N)\circ\partial_n\right)$$
and therefore
\begin{eqnarray*}S\otimes_RExt^n_\Lambda(M,N)&=&
S\otimes_R\left(\frac{Hom_\Lambda(\Omega^nM,N)}{Hom_\Lambda(P_{n-1},N)
\circ\partial_n}\right)\\
&=&
\frac{\left(S\otimes_RHom_\Lambda(\Omega^nM,N)\right)}{
\left(S\otimes_R\left(Hom_\Lambda(P_{n-1},N)\circ\partial_n\right)\right)}\\
&=&\frac{Hom_{S\otimes_R\Lambda}(S\otimes_R\Omega^nM,S\otimes_RN)}{
Hom_{S\otimes_R\Lambda}(S\otimes_RP_{n-1},S\otimes_RN)\circ (1_S\otimes\partial_n)}\\
&=&Ext^n_{S\otimes_R\Lambda}(S\otimes_RM,S\otimes_RN)
\end{eqnarray*}
for all $n\in\N$, natural in $M$ and $N$.

\medskip

The case $k=N+1$ shows then
$$(C_N(X_S)[N+1],C_N(Y_S)[N+1])_S=S\otimes_R(C_N(X)[N+1],C_N(Y)[N+1])$$
and
$$(C_N(X_S)[N],C_N(Y_S)[N+1])_S=S\otimes_R(C_N(X)[N],C_N(Y)[N+1]).$$
By the case for bounded complex of projectives we get that
the natural morphism is an isomorphism for
$$(\tau_NX_S[1],C_N(Y_S)[N+1])_S\simeq S\otimes_R(\tau_NX[1],C_N(Y)[N+1])$$
and
$$(\tau_NX_S,C_N(Y_S)[N+1])_S\simeq S\otimes_R(\tau_NX,C_N(Y)[N+1]).$$
Therefore also
$$(X_S,C_N(Y_S)[N+1])_S\simeq S\otimes_R(X,C_N(Y)[N+1])$$
and by the very same arguments, using $k=N$, we get
$$(X_S,C_N(Y_S)[N])_S\simeq S\otimes_R(X,C_N(Y)[N]).$$
This shows that we get isomorphisms
in the two left and the two right vertical morphisms of $(\dagger)$
and hence also the central
vertical morphism is an isomorphism.
Hence $$(X_S,Y_S)_S\simeq S\otimes_R(X,Y)$$
and the lemma is proved.
\dickebox

\begin{Theorem}\label{abstractisos}
Let $R$ be a commutative Noetherian ring, let $S$ be a
commutative Noetherian $R$-algebra
and suppose that $S$ is a faithfully flat $R$-module.
Suppose $S\otimes_Rrad(R)=rad(S)$.
Let $\Lambda$ be a Noetherian $R$-algebra,
let $X$ and $Y$ be two objects of
of $D^b(\Lambda)$ and suppose that
$End_{D^b(\Lambda)}(X)$ is a
finitely generated $R$-module.
Then
$$S\otimes_R^{\mathbb L}X\simeq S\otimes_R^{\mathbb L}Y\Leftrightarrow X\simeq Y.$$
\end{Theorem}

\begin{Rem}
We observe that if $R$ is local and $S=\hat R$ is the
$rad(R)$-adic completion, then  $S$ is faithfully flat as
$R$-module and $S\otimes_Rrad(R)=rad(S)$.
\end{Rem}

Proof of Theorem~\ref{abstractisos}.
According to the hypotheses we now suppose that
$End_{D^b(\Lambda)}(X)$ and $End_{D^-(\Lambda)}(Y)$ are
finitely generated $R$-module and that $S\otimes_Rrad(R)=rad(S)$.
Since $S$ is flat over $R$, tensor product of $S$ over $R$
is exact and we may replace the left derived tensor product
by the ordinary tensor product.
We only need  to show ''$\Rightarrow$'' and assume therefore that
$X$ and $Y$ are in $K^{-,b}(\Lambda-proj)$, and that
$S\otimes_RX$ and $S\otimes_RY$ are isomorphic.

\medskip

Let $X_S:=S\otimes_RX$ and  $S\otimes_RY=:Y_S$  in $D^b(S\otimes_R\Lambda)$
to shorten the notation and denote by $\varphi_S$ the isomorphism $X_S\lra Y_S$.
Since then $X_S$ is a direct factor of $Y_S$ by means of $\varphi_S$, the mapping
$$\varphi_S=\sum_{i=1}^ns_i\otimes \varphi_i:X_S\lra Y_S$$
for $s_i\in S$ and $\varphi_i\in Hom_{D^b(\Lambda)}(X,Y)$ has a left inverse
$\psi:Y_S\lra X_S$ so that $$\psi\circ\varphi_S=id_{X_S}.$$

Then, $$0\lra rad(R)\lra R\lra R/rad(R)\lra 0$$
is exact and since $S$ is flat over $R$ we get that
$$0\lra S\otimes_Rrad(R)\lra S\lra S\otimes_R(R/rad(R))\lra 0$$
is exact. This shows
$$S\otimes_R(R/rad(R))\simeq S/(S\otimes_Rrad(R)).$$
By hypothesis we have $S\otimes_Rrad(R)=rad(S)$, identifying
canonically $S\otimes_RR\simeq S$.
Then there are $r_i\in R$ so that $1_S\otimes r_i-s_i\in rad(S)$
for all $i\in\{1,\dots,n\}$.

Put $$\varphi:=\sum_{i=1}^nr_i\varphi_i\in Hom_{D^b(\Lambda)}(X,Y).$$
Then
\begin{eqnarray*}
\sum_{i=1}^n\psi\circ(1_S\otimes (r_i\varphi_i))-1_S\otimes id_X&=&
\sum_{i=1}^n
\left(\psi\circ(1_S\otimes r_i\varphi_i)-
\psi\circ(s_i\otimes\varphi_i)\right)\\
&=&\sum_{i=1}^n
\left(1_S\otimes r_i-s_i\right)\cdot \left(\psi\circ(id_S\otimes\varphi_i)\right)\\
&\in &\left(rad(S)\otimes_REnd_{D^b(\Lambda)}(X)\right)
\end{eqnarray*}
and since $End_{D^b(\Lambda)}(X)$ is a Noetherian $R$-module,
using Nakayama's lemma we obtain that
$\psi\circ(\sum_{i=1}^n1_S\otimes r_i\varphi_i)$ is invertible
in $S\otimes_REnd_{D^b(\Lambda)}(X)$. Hence
$id_S\otimes_R\varphi$ is left split and therefore
$$
X_S\stackrel{id_S\otimes_R\varphi}{\lra}Y_S\lra
C(id_S\otimes_R\varphi)\stackrel{0}{\lra}X_S[1]
$$
is a distinguished triangle, with $C(id_S\otimes_R\varphi)$
being the cone of $id_S\otimes_R\varphi$.
However, $$C(id_S\otimes_R\varphi)=S\otimes_RC(\varphi)$$
and hence
$$
X_S\stackrel{id_S\otimes_R\varphi}{\lra}Y_S\lra
S\otimes_RC(\varphi)\stackrel{0}{\lra}X_S[1]
$$
is a distinguished triangle.

Since $\varphi_S$ is an isomorphism, $\varphi_S$ has a
right inverse $\chi:Y_S\lra X_S$ as well.
Now, since $X_S\simeq Y_S$, since $S$ is faithfully flat over $R$,
and since $End_{D^b(\Lambda)}(X)$ is finitely generated as $R$-module,
using Lemma~\ref{stoensorhomequalhomtensors} we obtain that
$End_{D^b(\Lambda)}(Y)$ is finitely generated as $R$-module as well.
The same argument as for the left inverse $\psi$ shows that
$(id_S\otimes\varphi)\circ\chi$ is invertible
in $S\otimes_REnd_{D^b(\Lambda)}(Y)$. Hence
$$
X_S\stackrel{id_S\otimes_R\varphi}{\lra}Y_S\stackrel{0}{\lra}
S\otimes_RC(\varphi)\stackrel{0}{\lra}X_S[1]
$$
is a distinguished triangle. This shows that $S\otimes_RC(\varphi)$
is acyclic, and hence $$0=H(S\otimes_RC(\varphi))=S\otimes_RH(C(\varphi)).$$
Since $S$ is faithfully flat over $R$ also $H(C(\varphi))=0$,
which implies that $C(\varphi)$ is acyclic
and therefore $\varphi$ is an isomorphism.

This proves the theorem.
\dickebox

\bigskip

Let $A$ be an algebra over a complete discrete valuation ring $R$ which is
finitely generated as a module over $R$.
We shall need a Krull-Schmidt theorem for the derived category of bounded
complexes over $A$.  This fact seems to be well-known, but
for the convenience of the reader we give a proof.

\begin{Prop}\label{KrullSchmidtforboundedprojectives}
Let $R$ be a complete discrete valuation ring and let $A$ be an
$R$-algebra, finitely generated as $R$-module.
Then the Krull-Schmidt theorem holds for $K^{-,b}(A-proj)$.
\end{Prop}

Proof.
We first show a Fitting lemma for $K^{-,b}(A-proj)$.

Let $X$ be a complex in $K^{-,b}(A-proj)$ and
let $u$ be an endomorphism of the complex $X$. Then $X=X'\oplus X''$
as graded modules, by Fitting's lemma in the version for
algebras over complete discrete valuation rings \cite[Lemma 1.9.2]{Benson}.
The restriction of $u$ on
$X'$ is an automorphism in each degree and the restriction of
$u$ on $X''$ is nilpotent modulo $rad(R)^m$ for each $m$.
Therefore $u$ is a diagonal matrix
$\left(\begin{array}{cc}\iota&0\\0&\nu\end{array}\right)$ in each degree
where $\iota:X'\lra X'$ is invertible, and $\nu:X''\lra X''$ is
nilpotent modulo $rad(R)^m$ for each $m$
in each degree.
The differential $\partial$ on $X$ is given by
$\left(\begin{array}{cc}\partial_1&\partial_2\\ \partial_3&\partial_4\end{array}\right)$
and the fact that $u$ commutes with $\partial$ shows that
$\partial_3\iota=\nu\partial_3$ and $\partial_2\nu=\iota\partial_2$.
Therefore, $\partial_3\iota^s=\nu^s\partial_3$ and
$\partial_2\nu^s=\iota^s\partial_2$ for all $s$. Since
$\nu$ is nilpotent modulo $rad(R)^m$ for each $m$ in each degree,
and $\iota$ is invertible, $\partial_2=\partial_3=0$.
Hence the differential of $X$ restricts to a differential
on $X'$ and a differential on $X''$.
Moreover, $X'$ and $X''$ are both
projective modules, since $X$ is projective.

Now, $X$, and therefore also
$X''$ is exact in degrees higher than $N$, say. We fix $m\in\N$
and obtain therefore that $u$ is nilpotent modulo $rad(R)^m$
in each degree lower than $N$. Let $M_m$ be the nilpotency degree.
Then, since $X''$ is exact in degrees higher than $N$,
modulo $rad(R)^m$ the restriction of the endomorphism
$u^{M_m}$ to $X''$
is homotopy equivalent to $0$ in degrees higher than $N$. We get therefore
that the restriction of $u$ to $X''$ is actually nilpotent modulo
$rad(R)^m$ for each $m$.

Hence, the endomorphism ring of an indecomposable object is
local and the Krull-Schmidt theorem is an easy consequence by the
classical proof as in \cite{Lang} or in \cite{Benson}.

This shows the proposition. \dickebox

\begin{Rem}
If $R$ is a field and $A$ is a finite dimensional $R$-algebra,
then we would be able to argue more directly.
Indeed, $X'={\rm im}(u^N)$  and $X''=\ker(u^N)$  for large enough $N$.
Then it is obvious that $X'$ and $X''$ are both subcomplexes of $X$.
Observe that $R$ may be a field in Proposition~\ref{KrullSchmidtforboundedprojectives}.
\end{Rem}

For the next Corollary we follow closely \cite{Rog72}.

\begin{Cor}\label{finallyfaithfulprojectiveextensions}
Let $R$ be a commutative semilocal Noetherian ring, let $S$ be a commutative $R$-algebra
so that $\hat S:=\hat R\otimes_RS$ is a faithful projective $\hat R$-module
of finite type.  Let $\Lambda$ be a Noetherian $R$-algebra,
finitely generated as $R$-module, and
let $X$ and $Y$ be two objects
of $D^b(\Lambda)$ and suppose that $End_{D^b(\Lambda)}(X)$ and
$End_{D^b(\Lambda)}(Y)$ are finitely generated $R$-module. Then
$$S\otimes_R^{\mathbb L}X\simeq S\otimes_R^{\mathbb L}Y\Leftrightarrow X\simeq Y.$$
\end{Cor}

Proof. If $S\otimes_R^{\mathbb L}X\simeq S\otimes_R^{\mathbb L}Y$ in
$D^b(S\otimes_R\Lambda),$
we get $\hat S\otimes_R^{\mathbb L}X\simeq \hat S\otimes_R^{\mathbb L}Y$ in
$D^b(\hat S\otimes_R\Lambda)$.
Since $R$ is semilocal with maximal ideals $m_1,\dots,m_s$
we get $\hat R=\prod_{i=1}^s\hat R_{m_i}$ for
the completion $\hat R_{m_i}$ of $R$ at $m_i$.
Now, $\hat S$ is projective faithful of finite type, and
so there are $n_1,\dots,n_s$ with
$$\hat S\simeq \prod_{i=1}^s(\hat R_{m_i})^{n_i}$$ and therefore
$\hat S\otimes_R^{\mathbb L}X\simeq \hat S\otimes_R^{\mathbb L}Y$ \mbox{ implies }
$$\prod_{i=1}^s(\hat R_{m_i})^{n_i}\otimes_R^{\mathbb L}X\simeq \prod_{i=1}^s(\hat R_{m_i})^{n_i}\otimes_R^{\mathbb L}Y.$$
Hence
$$(\hat R_{m_i}\otimes_R^{\mathbb L}X)^{n_i}\simeq
(\hat R_{m_i}\otimes_R^{\mathbb L}Y)^{n_i}$$
for each $i$, and therefore by
Proposition~\ref{KrullSchmidtforboundedprojectives}
$$\hat R_{m_i}\otimes_R^{\mathbb L}X\simeq \hat R_{m_i}\otimes_R^{\mathbb L}Y$$
for each $i$. By Theorem~\ref{abstractisos} we obtain $X\simeq Y$.
\dickebox

\bigskip

We get cancellation of factors from this statement.

\begin{Cor}
Under the hypothesis of Theorem~\ref{abstractisos} or of
Corollary~\ref{finallyfaithfulprojectiveextensions}
we get $X\oplus U\simeq Y\oplus U$ in $D^b(\Lambda)$
implies $X\simeq Y$.
\end{Cor}

Proof. This is clear by Corollary~\ref{finallyfaithfulprojectiveextensions}
in combination with Proposition~\ref{KrullSchmidtforboundedprojectives}.
\dickebox

\begin{Rem}
In \cite{JSZ} we developed a theory to roughly speaking
parameterise geometrically objects in $D^b(A)$ by
orbits of a group action on a variety.
For this purpose we need to assume that $A$ is a finite dimensional algebra over
an algebraically closed field $K$, so that it is possible to
use arguments and constructions from algebraic geometry.
Using Theorem~\ref{abstractisos} we can extend the theory to non algebraically closed fields
$K$ as well.
\end{Rem}

\paragraph{\bf Acknowledgement:} I wish to thank the referee for
useful suggestions which lead to substantial improvements.

\end{document}